\documentclass[12pt, reqno]{amsart}
\usepackage{amsmath, amstext, amsbsy, amssymb, amscd}

\newtheorem{theorem}{Theorem}[section]
\newtheorem{lemma}[theorem]{Lemma}

\newtheorem{corollary}[theorem]{Corollary}
\theoremstyle{definition}

\theoremstyle{remark}
\newtheorem{remark}[theorem]{Remark}

\numberwithin{equation}{section}

\newcommand{\C}{ \mathbb C }

\newcommand{\Supp}{{\rm Supp}}

{\vskip-\lastskip\medskip
  \noindent
  {\em #1.}\enspace
  }%
{\qed\par\medskip
  }

\begin{document}

 \title[On the  existence of  nontrivial threefolds with vanishing Hodge cohomology]
        {On the existence of  nontrivial threefolds with vanishing Hodge cohomology}

 \author[Jing  Zhang]{Jing  Zhang}
\address{Department of Mathematics, University of 
Missouri, Columbia, MO
65211, USA}
\email{zhangj@math.missouri.edu}

\begin{abstract}
We  analyze the structure of  
  the algebraic manifolds $Y$ of dimension 3 with 
$H^i(Y, \Omega^j_Y)=0$
for all $j\geq 0$, $i>0$ and 
$h^0(Y, {\mathcal{O}}_Y) > 1$, 
by showing the deformation invariant of  some
 open surfaces.  Secondly, we   show   when
  a smooth threefold  with nonconstant regular functions
 satisfies the  vanishing Hodge cohomology. 
As an application, we  
 prove the existence  of nonaffine and nonproduct threefolds  
 $Y$
with  this property 
by constructing a family of a certain type of
open surfaces parametrized by the affine curve $\C-\{0\}$
 such that the corresponding  
smooth completion
$X$ has  Kodaira dimension 
$-\infty$ and $D$-dimension 1, where $D$ is the effective boundary divisor
with support $X-Y$.  
\end{abstract}

\maketitle
\date{}
%\tableofcontents
%%
%%
%%
%%
%%
%%
%%
\section{Introduction}

We  study  the structure of  algebraic  manifolds  $Y$  of dimension 3
with     $H^i(Y, \Omega^j_Y)=0$
for all $j\geq 0$, $i>0$.  Originally this question was 
raised by J.-P. Serre for complex manifolds \cite{Se}. Since by Serre duality
$Y$ is not complete,   
$Y$ is affine if  it  is a curve (\cite{H2}, page 68).
If  $Y$ is a surface, it was classified by Mohan Kumar \cite{Ku} 
(see the following theorem in this section). We are interested in three dimensional 
case. 
Suppose that $X$ is a smooth 
completion of $Y$. If there are nonconstant regular functions on $Y$,
i.e., $h^0(Y, {\mathcal{O}}_Y) > 1$, 
then  $Y$ contains no complete curves and 
 the boundary    is connected \cite{Zh}. 
Therefore we may assume that the boundary is of pure codimension 1
by suitable blowing up subvarieties on the boundary.
Let $D$ be an effective divisor with simple normal crossings 
(\cite{KoM}, page 5) such that
   supp$D=X-Y$. The condition
$h^0(Y, {\mathcal{O}}_Y) > 1$ is equvalent to $\kappa(D, X)>0$. Here we use the 
standard definition of $D$-dimension due to Iitaka.    
If for all integers  $m> 0$ we have
     $H^0(X, {\mathcal{O}}_X(mD))=0$, then we define 
     the  $D$-dimension of $X$, denoted by $\kappa (D, X)$, to be $-\infty$.
     If  $h^0(X, {\mathcal{O}}_X(mD))\geq 1$ for some $m$, 
     choose a basis $\{f_0, f_1, \cdot \cdot\cdot, f_n\}$
     of the linear space 
     $H^0(X, {\mathcal{O}}_X(mD))$, it defines a rational 
     map 
     $\Phi _{mD}$
     from $X$ to the projective space 
     ${\Bbb{P}}^n$ by sending a point $x$ on $X$ to
     $(f_0(x), f_1(x), \cdot \cdot\cdot, f_n(x))$ in ${\Bbb{P}}^n$.  
     Then we define
     $\kappa (D, X)$ to be the maximal dimension of the images
      of the rational map  $\Phi _{mD}$, i.e., 
      $$ \kappa (D, X)= \max_m\{\dim (\Phi _{mD}(X))\}. 
      $$
Let $K_X$ be the canonical divisor of $X$, then the Kodaira 
    dimension of $X$ is  the $K_X$-dimension of $X$, denoted by 
    $\kappa(X)$, i.e.,  
    $$\kappa(X)=\kappa(K_X, X). 
    $$  

Before we state our theorems, we need Mohan Kumar's result for surfaces.

\noindent  
{\bf{Theorem (Mohan Kumar)}}
 {\it  Let
$Y$ be a smooth algebraic surface over $\Bbb{C}$
with   $H^i(Y, \Omega^j_Y)=0$   for all $j\geq 0$ and $i>0$,
then  $Y$ is one of the following 
 
      (1) $Y$ is affine.

      (2) Let $C$ be an elliptic curve and $E$ the unique nonsplit 
      extension of $\mathcal{O}$$_C$ by itself.  
      Let ${X=\Bbb{P}}_C(E)$ and  $D$ be the canonical section, then $Y=X-D$.

       (3) Let $X$ be a projective rational surface with an effective 
       divisor $D=-K$ with $D^2=0$, $\mathcal{O}$$(D)|_D$ be nontorsion and 
       the dual graph of $D$ be $\tilde{D}_8$ or $\tilde{E}_8$, then $Y=X-D$.}

\begin{theorem} If $Y$ is an algebraic manifold of dimension 3 with  
$H^i(Y, \Omega^j_Y)=0$ for all $j\geq 0$ and $i>0$
and $h^0(Y, {\mathcal{O}}_Y) > 1$, then  we have a surjective 
morphism  from  $Y$ to a smooth affine curve $C$ 
such that all smooth fibres are of the same 
type, i.e., exactly one of the three types of open surfaces in Mohan Kumar's
classification.  
Moreover, if one fibre is not affine, then 
$X$ has  Kodaira dimension 
$-\infty$ and $D$-dimension 1. 
\end{theorem} 

It is well-known that the type 
(2) and type (3) projective surfaces are rigid. However, the rigidity of
the projective surfaces does not imply the rigidity of the open
surfaces. The problem  is that if  a surface is affine, then 
its smooth completion 
can be any projective surface. In particular, it can be
 type (2) or type (3)
projective surface. 
More precisely, assume that a type (3) projective surface $X_0$
deforms to a projective surface $X_1$, then $X_0$ and $X_1$ have 
the same minimal model
\cite{I4}, \cite{BaPV}, Chapter VI, Theorem 8.1. 
Let $S_0$ and $S_1$ be the corresponding open surfaces in $Y$ 
contained in
$X_0$ and $X_1$ respectively, then we have 
 $H^i(S_0, \Omega^j_{S_0})=H^i(S_1, \Omega^j_{S_1})=0$ 
for all $j\geq 0$ and $i>0$ \cite{Zh}. Since  both affine surface
and type (3) open surface satisfy this condition, even 
though $X_0$ is isomorphic to $X_1$, in priory, 
$S_0$ and $S_1$ may not be of the same type.  
So we need to rule out 
the following case: some isolate fibre is affine but general fibres are not 
affine. We will carefully analyze 
how the cohomology of the sheaves 
${\mathcal{O}}_X(nD)$ changes
when restricted to each fibre
to obtain the  deformation invariant of the open surfaces.

\begin{theorem}  With the same assumption as in the above theorem,
  if one smooth fibre $S_0$
of $f|_Y$ over $t_0\in C$ is affine, then 
by removing finitely many fibres  $S_1$, $S_2$,$\cdot, \cdot, \cdot$,
 $S_m$ from Y,
 the new threefold
$Y'=Y- \cup S_i$
is affine. 
\end{theorem}

When restricted to a fibre,  if the global divisor $D$ on $X$ is ample, then
the above theorem is trivial by \cite{KoM},  Proposition 1.41. However,
if an open fibre is affine, we only know that its boundary on the 
corresponding projective 
surface is the support of an ample divisor on the surface. There is 
 no guarantee that
this ample divisor on the fibre can be extended to a global divisor
on $X$.  We will use Goodman and Hartshorne's result (Lemma 3.1) to
transfer the cohomology condition on the open fibre to the closed
fibre in order to apply upper semicontinuity theorem.

  Let  $\bar{C}$   
be a smooth projective curve containing $C$.
Let $F_n=\Omega^j_X\otimes {\mathcal{O}}_X(nD)$. 
Now we do not assume  $H^i(Y, \Omega^j_Y)=0$ for 
all $j\geq 0$ and $i>0$. We want to know whether 
$Y$ satisfies this condition if every fibre $S$
satisfies it, i.e., $H^i(S, \Omega^j_S)=0$ for all $j\geq 0$ and $i>0$. 
We know that if globally $Y$ is such a threefold,
then each fibre must satisfy the same vanishing condition \cite{Zh}.
The converse is   very  subtle. Assume that each fibre and the base 
satisfy some property in a fibre space, then  globally   the property may fail. 
 A famous example is Skoda's counterexample  \cite{Sk}
for Serre's question \cite{Se}: Is the total space of a holomorphic 
    fibre bundle with Stein base $Z$ and Stein fibre $F$ a Stein manifold?  
In order to prove that the vanishing Hodge cohomology holds for  $Y$,
we  will  first prove the locall freeness of the higher direct images 
$R^if_*F_n$ for $n\gg 0$. The local freeness is interesting 
on its own.

\begin{theorem} If we have the  commutative diagram
\[
  \begin{array}{ccc}
    Y                           &
     {\hookrightarrow} &
    X                                 \\
    \Big\downarrow\vcenter{%
        \rlap{$\scriptstyle{f|_Y}$}}              &  &
    \Big\downarrow\vcenter{%
       \rlap{$\scriptstyle{f}$}}      \\
C        & \hookrightarrow &
\bar{C}
\end{array}
\]
such that f is proper and surjective and  
 each fibre $X_t$ over $t\in C$ is of type (2)
projective surface, then $R^if_*F_n|_C$ is locally free for all
$i\geq 0$ and
$n\gg 0$. Therefore   $H^i(Y, \Omega^j_Y)=0$ for all $j\geq 0$ 
    and $i>0$. 
\end{theorem} 

If we  also  assume that every horizontal divisor 
$D_i$ (i.e., $f(D_i)=\bar{C}$)   intersects each 
smooth fibre $X_t=f^{-1}(t)$ 
over $t\in C$ 
with one  
prime divisor  on $X_t$, then for type 
(3) fibres, the theorem still holds. We add this technical 
condition because a prime component of $D$ might intersect 
some fibre with two or more curves.

\begin{theorem} In the above commutative diagram,
if each fibre $X_t$ over $t\in C$ is of type (3)
projective surface, then $R^if_*F_n|_C$ is locally free for all $i\geq 0$
and $n\gg 0$.  Furthermore,   $H^i(Y, \Omega^j_Y)=0$ for all $j\geq 0$ 
    and $i>0$. 
\end{theorem} 

\begin{corollary} If we have a surjective morphism from
a smooth threefold $Y$ to a smooth affine curve $C$ such
that every fibre is smooth and the above diagram commutes, then 
$H^i(Y, \Omega^j_Y)=0$ for all $j\geq 0$ 
    and $i>0$ if and only if 
$H^i(S, \Omega^j_S)=0$ for every fibre $S$ and  all $j\geq 0$ 
    and $i>0$.
\end{corollary}

As a consequence of  above theorems, we have the following existence theorem.
\begin{theorem} There exist nonaffine and nonproduct threefolds $Y$ with 
$H^i(Y, \Omega^j_Y)=0$ for all $j\geq 0$ 
    and $i>0$.
\end{theorem}

We will prove  these  theorems in the following sections. 
The proof of Theorem 1.3 and Theorem 1.4 is similar. We will just 
prove Theorem 1.3. 
\\

\noindent 
{\bf{Question}}\quad  
Are the threefolds $Y$  Stein in Theorem 1.3 and Theorem 1.4?\\

\noindent 
{\bf  Convention}\quad Unless otherwise explicitly mentioned, 
we always use Zariski topology, i.e.,  
     an open set means a Zariski open set.\\

\noindent 
{\bf{Acknowledgments}}  \quad  I would like to  express my  thanks to 
   the following 
      professors for  helpful discussions: Michael Artin,  Steven Dale Cutkosky,
     Dan Edidin,   N.Mohan  Kumar, 
      Zhenbo Qin, A. Prabhakar  Rao,  David Wright and Qi Zhang.

%\tableofcontents
%%
%%
%%
%%
%%
%%
%%
\section{Proof of  Theorem 1.1}

\noindent 
{\bf{Theorem}{(Iitaka)}} 
{\it Let $X$ be a normal projective variety and let 
$D$ be an effective divisor on $X$. Then  there exist two positive
numbers $\alpha$ and $\beta$ such that for all sufficiently 
large $n$ we have
$$  \alpha n^{\kappa(D, X)}
\leq h^0(X, {\mathcal{O}}_X(nD))
\leq \beta n^{\kappa(D, X)}. 
$$ }

For the proof of  Iitaka's theorem, see Lecture 3 \cite{I1}
or Theorem 8.1 \cite{U}.

The following  two 
lemmas  are known \cite{Ku}. 

\begin{lemma} Let $S$ be a smooth  open surface
with $H^i(S, \Omega^j_S)=0$
for all $j\geq 0$ and $i>0$. Let  $\bar{S}$
be  a  smooth projective surface  containing $S$
and $G$ be the divisor in Mohan Kumar's Theorem.
Then there are three cases.

(1)  If $S$ is affine, then   $\kappa(G, \bar{S})=2$ and 
$$h^0(\bar{S}, {\mathcal{O}}_{\bar{S}}(nG)) \geq cn^2$$
for some positive constant $c$ and $n\gg 0$.

(2) If $S$ is of type (2), then  
$$h^0(\bar{S}, {\mathcal{O}}_{\bar{S}}(nG))=
h^1(\bar{S}, {\mathcal{O}}_{\bar{S}}(nG))=1, \quad 
h^2(\bar{S}, {\mathcal{O}}_{\bar{S}}(nG))=0$$

for all $n\gg 0$. 

(3) If $S$ is of type (3), then 
$$h^0(\bar{S}, {\mathcal{O}}_{\bar{S}}(nG))=1, \quad 
h^1(\bar{S}, {\mathcal{O}}_{\bar{S}}(nG))=
h^2(\bar{S}, {\mathcal{O}}_{\bar{S}}(nG))=0$$
for all $n\gg 0$. 
\end{lemma}

$Proof.$ (1) Since $S$ is affine, 
$\bar{S}-S$ is support of an ample divisor
$A$ by Goodman's theorem \cite{H2}, page 69. So 
$\kappa(A, \bar{S})=\kappa(G, \bar{S})=2$ \cite{I6} and 
\cite{Ba}, Chapter 14.
The  estimate  is obvious by Iitaka's Theorem.

(2) The equalities follow from 
Lemma 1.8, \cite{Ku} and the following Lemma 2.2(1).

(3) See Lemma 1.8,  3.1 \cite{Ku} and the following Lemma 2.2(2).
\begin{flushright}
 Q.E.D. 
\end{flushright}

\begin{lemma}  With the preceding notation,  we have

(1)  If $\bar{S}$ is of type (2), then $G^2=0$, $K_{\bar{S}}=-2G$,
 $p_g=0$ and $q=1$.

(2) If $\bar{S}$ is of type (3), then $G^2=0$, $K_{\bar{S}}=-G$,
 $p_g=q=0$.
\end{lemma}

$Proof.$ (1) This is a standard result for the ruled surface 
over an elliptic curve. The proof can be found in \cite{H1},
 Chapter V, Section 2
or \cite{Ku}. 

(2) See Lemma 1.6  and 3.1   \cite{Ku}.
\begin{flushright}
 Q.E.D. 
\end{flushright}

Let $f$:  $X\rightarrow Z$ be a morphism between varieties (schemes)
with $Z$ connected.  Let   $z_0\in Z$, $k(z_0)=K$, and 
$X_{z_0}\cong X_0$. Then the other fibres $X_z$ of $f$ are called 
deformations of $X_0$, [H1], page 89. In the  proof of Theorem 1.1,
 the deformation
of a nonsingular  complex surface $X_0$
means the following by the same notation:
 Both $X$ and $Z$ are smooth
and $f$ is  surjective, proper and  flat morphism 
(i.e.,  ${\mathcal{O}}_{X, x}$ is a flat 
${\mathcal{O}}_{Z, f(x)}$-module for all $x\in X$)
such that
the fibre over $z_0\in Z$, $X_{z_0}\cong X_0$, \cite{BaPV}, page 36. 
By [I4], we know  that
the deformation of a rational surface is again rational. 
By Theorem (8.1), Chapter VI, [BaPV],  
the deformation of  a  ruled surface over a smooth curve 
of genus $g\geq 1$
is also of the same type, i.e., has the same minimal model.

We need Kodaira's  stability of  $(-1)$-curves. It is Theorem 5 in [Kod1].

\noindent 
{\bf{Theorem}{(Kodaira)}} 
{\it 
Let $f$:  $X\rightarrow Z$ be a surjective, proper holomorphic
map which is flat. If for some point 
$0\in Z$ the fibre $X_0$ contains a  $(-1)$-curve $E_0$, then 
there is an open neighborhood  (in complex topology) $U$ of 0
in $Z$ and a closed and connected submanifold $E$ of $f^{-1}(U)$
such that  $E\cap X_0=E_0$ and such that   
$E\cap X_t=E_t$ is a  $(-1)$-curve  for every  $t\in U$. }

Further, in Kodaira's Theorem,    
there is  a  $g$: $X' \rightarrow U$,  a surjective, flat,  
proper holomorphic
map such that the following diagram commutes: 
\[
  \begin{array}{ccc}
    X                           &
     {\stackrel{h}{\longrightarrow}} &
    X'                                 \\
    \Big\downarrow\vcenter{%
        \rlap{$\scriptstyle{f}$}}              &  &
    \Big\downarrow\vcenter{%
       \rlap{$\scriptstyle{g}$}}      \\
U        &{\stackrel{\approx} {\longrightarrow} }&
U
\end{array}
\]
where  $h|_{X_t}: X_t\rightarrow X_t'$ is the blowing down of $E_t$. 
Let us state the contraction part precisely. The proof
 is due to Suwa [I5], Appendix 1.

\noindent 
{\bf{Theorem}{(Suwa)}} 
{\it  Let $X$ and $Z$ be complex manifolds, 
and let $f$
be a proper, surjective and flat holomorphic map from $X$ to $Z$,
such that every fibre $X_z$ is a smooth surface. If there exists
a complex submanifold $E$ of $X$ such that its restriction to
$X_z$: $E_z=E\cap X_z$ is an irreducible exceptional curve of the first 
kind on $X_z$ at any $z\in Z$, then we can construct a complex manifold $X'$,
which is proper over $Z$, and a holomorphic map $h$: $X\rightarrow X'$
over $Z$, such that $h|_{X_z}$: $X_z\rightarrow X'_z$  shrinks  $E_z$
to a point in  $X_z'$ for every point $z\in Z$, and such that 
$h|_{X-E}$:  $X-E\rightarrow X'-h(E)$
is biholomorphic. }

\noindent
{\bf{Upper Semicontinuity Theorem (Grauert, Grothendieck) }} {\it
Let $f$: $X\rightarrow Z$  be a proper morphism of noetherian 
schemes, ${\mathcal{F}}$ a coherent sheaf on $X$, flat over $Z$,
then

(1) The  $i$-th    direct image 
 $R^if_*{\mathcal{F}}$  is  a coherent sheaf on $Z$ for any
 nonnegative integer $i$.
 
 (2) Let  ${\mathcal{F}}_z={\mathcal{F}}|_{X_z}$, i.e., 
 the sheaf ${\mathcal{F}}$ restricted to the fibre 
 $X_z=f^{-1}(z)$, then the function 
 $$ d_i(z)= h^i(X_z,  {\mathcal{F}}_z)=
 \dim_{k(z)}H^i(X_z, {\mathcal{F}}_z)
 $$
 is upper semi-continuous on $z$. That is, for any 
 $n\in  {\Bbb{Z}}$, the set 
 $\{z\in Z: d_i(z) \geq n\}$  is a closed set, where 
 $k(z)={\mathcal{O}}_z/{\mathcal{M}}_z$, the residue field
 at the point $z$.  

 (3) The Euler characteristic of the restriction sheaf 
 ${\mathcal{F}}_z$ 
 $$ \chi ({\mathcal{F}}_z)= 
 \sum (-1)^i {\dim}_{k{(z)}}H^i(X_z, {\mathcal{F}}_z)
 $$
 is locally constant on $Z$. 
 
 (4) The following statements are equivalent:

 \quad ({$\mbox{i}$}) $h^i(X_z,  {\mathcal{F}}_z)$
 is a constant function on $Z$,
 
 \quad ({$\mbox{ii}$}) $R^if_*{\mathcal{F}}$  is  locally free 
  sheaf on $Z$, and for all $z\in  Z$, the natural map 
  $$ R^if_*{\mathcal{F}}\otimes _{{\mathcal{O}}_z} k(z)
  \longrightarrow 
  H^i(X_z, {\mathcal{F}}_z)
  $$ 
  is an isomorphism. 
  
  In addition, if these conditions are satisfied,
  then
   $$ R^{i-1}f_*{\mathcal{F}}\otimes _{{\mathcal{O}}_z} k(z)
  \longrightarrow 
  H^{i-1}(X_z, {\mathcal{F}}_z)
  $$ 
  is an isomorphism  for all $z\in  Z$. 
 }

For a proof, see   \cite{Mu}, page 46-53.

In this section, from now on,  we assume that the 
condition of Theorem 1.1 holds. 
 Theorem 1.1 is a direct consequence of  the following Lemmas.

\begin{lemma}If one smooth fibre
$S_{t_0}=S_0$ is of type (2) or (3) 
open surface
in Mohan Kumar's classification,
then there is an affine open set $U$ such that 
 $S_t=f^{-1}(t)-D$ over every $t\in U$ is of the same type.
\end{lemma}

$Proof.$  Consider the commutative  diagram  in \cite{Zh}
\[
  \begin{array}{ccc}
    Y                           &
     {\hookrightarrow} &
    X                                 \\
    \Big\downarrow\vcenter{%
        \rlap{$\scriptstyle{f|_Y}$}}              &  &
    \Big\downarrow\vcenter{%
       \rlap{$\scriptstyle{f}$}}      \\
C        & \hookrightarrow &
\bar{C},
\end{array}
\]
where $f$ is proper and surjective and for all $t\in C$, 
$X_t=f^{-1}(t)$ is a smooth projective surface. 
The minimal model of $X_t$ is the same as the minimal model 
 of type (2) or type (3) 
surface but may contains exceptional curves of the 
first kind.

Note that 
$S_0$ is not affine.  Let $X_0=f^{-1}(t_0)$.
 By Lemma 2.1  and \cite{U} Lemma 5.3,
even though the divisor $D_0=D|_{X_0}$ contains exceptional curves 
of the  first kind,  we still have 
$H^0(X_0, {\mathcal{O}}(nD_0))=\Bbb{C}$ for all nonnegative integer $n$. 
Let  $D_t=D|_{X_t}$, then $D_t$ is a connected curve on $X_t=f^{-1}(t)$ since 
$X_t$ is smooth  and  $H^i(S_t, \Omega^j_{S_t})=0$ \cite{Ku}Lemma 1.4. 
By upper semi-continuity, there is an affine open set $U$ in $C$ such that
$H^0(X_t, {\mathcal{O}}(nD_t))
=\Bbb{C}$ since every $D_t$ is effective. 
Therefore  
every fibre $S_t$ over $t\in U$ is not affine by Lemma 2.1. 

Secondly, if  $S_0$ is of type (2) open surface 
 in Mohan Kumar's classification, then  
$p_g(X_0)=h^2({\mathcal{O}}_{X_0})=0$ and 
$q(X_0)=h^1({\mathcal{O}}_{X_0})=1$  
(Lemma  2.2). Here the boundary divisor 
$D_0=D|_{X_0}$ may contain exceptional curves so $X_0$
may not be minimal. But $p_g$ is
birational invariant and $q$ is bimeromorphic invariant \cite{BaPV}, page 107
and \cite{H1}, page 181. 
Since $R^if_*({\mathcal{O}}_X)$
and   $R^if_*({\mathcal{O}}_X(K_X))$ are locally free for all
$i\geq 0$ \cite{Kol1,  Kol2}, again by upper semi-continuity, 
$p_g(X_t)=0$ and  $q(X_t)=1$ for every $t\in C$. Now $X_0$ has the  minimal 
model of a
 ruled surface over an elliptic curve, by the  classification theorem (1.1),
page 243
and deformation theorem (8.1), page 263, Chapter VI [BaPV], 
there is an affine   open set $U$ such that for
every $t\in U$, 
$X_t$ has the same minimal 
model as $X_0$
  in Mohan Kumar's theorem, i.e.,
type (2) projective surface.

 Similarly, if $S_0$ is of 
the third type, then there is an affine open set $U$ such that 
 $S_t$ over every $t\in U$ is of the same type 
 since  the deformation  of a rational surface is still
  rational \cite{I4}. 
\begin{flushright}
 Q.E.D. 
\end{flushright}

\begin{remark} If $S$ is of type (2) open surface in Mohan Kumar's Theorem,
then any point on $S$ cannot be contained in any exceptional 
curve of $\bar{S}$, where $\bar{S}$  is a smooth completion of $S$. 
So if $\bar{S}$ is not minimal, then all exceptional curves 
are contained in the boundary $\bar{S}-S$. 
\end{remark}

\begin{lemma}If there is an affine open set $U$ in $C$ such that 
for every $t\in U$, $t\neq t_0$, $S_t=f^{-1}(t)-D$ is of type (2)
open surface
 of Mohan Kumar's classification, 
where $t_0$ is a fixed point of $U$, then $S_0$ must be of the same
type surface. 
\end{lemma}
$Proof.$  First, $S_0$ cannot be of type (3) since $X_t=f^{-1}(t)$
is not rational and the deformation of rational surfaces is still 
rational \cite{I4}. 
We know that 
there are three possible smooth fibres \cite{Ku, Zh}. So  
we only need to prove that $S_0$ is not affine.  
It suffices to prove that $h^0(X_0, {\mathcal{O}}_{X_0}(nD_0))$
is bounded for all $n$. In fact, in our case, it is 1. 
Here $X_0=f^{-1}(t_0)$, $D_0=D|_{X_0}$ and $S_0=X_0-D_0$.

By Kodaira's Stability Theorem of $(-1)$-curves and 
Suwa's Theorem, we may assume that $D_0$  has no exceptional 
curve of the first kind. So there is a small open 
set $V$ in $C$ (complex topology), for all points
$t\in V$, $D_t=D|_{X_t}$ has no exceptional curves of the first kind.
In fact, if  there is $t_1\in V$, $D_{t_1}$
has a component $E_1$, such that $E_1$ is an exceptional curve 
of the first kind, then $E_1^2=E_1\cdot K_{X_1}=-1$, where 
$X_1=f^{-1}(t_1)$. There is a prime  component  $G$
of $D$ in $X$ such that $E_1\subset G$. Let $E_t=G|_{X_t}$ for $t\in V$,
then    by upper semi-continuity,
  the Euler characteristic  of 
     ${\mathcal{O}}_{X_t}(nE_t)$ 
     is constant for every $t\in V$ and every $n\geq 0$.
So for any $n\geq 0$,  we  have 
$$ \chi ({\mathcal{O}}_{X_t}(nE_t))
=\chi ({\mathcal{O}}_{X_0}(nE_0)).
$$
 By Riemann-Roch formula, for all $n\geq 0$, we have
$$  \frac{1}{2}n^2E_t^2-\frac{1}{2}nE_t\cdot K_{X_t}
=\frac{1}{2}n^2E_1^2-\frac{1}{2}nE_1\cdot K_{X_1}.
$$ 
So  $E_t^2=E_t\cdot K_{X_t}=-1$ for all $t\in V$. In particular,
$E_0^2=E_0\cdot K_0=-1$. This is impossible since $D_0$ has no $(-1)$-curves
by our assumption. Thus for all $t\in V$ and $t\neq t_0$,
$X_t$ is type (2) surface, i.e., a minimal ruled surface over an 
elliptic curve.  But $D_t$ may not be a prime divisor.
 Let $D_t'$ be the elliptic 
curve (a section) as in Mohan Kumar's classification,
 then there is a positive integer 
$n(t)$, depending on $t$ such that $D_t=n(t)D_t'$. Since the function
$n(t)$ is discrete, there is a dense  subset
$B$ in $V$ such that $n(t)$ is a constant $c$ for all $t\in B$. Let  
$t_1\in V-B$ and $K_1=K_{X_1}$.   
Considering the  divisor $D+cK_X$, when  restricted to the 
corresponding fibre
 $X_1=f^{-1}(t_1)$,
  by upper  semi-continuity, we have
  $$  h^0(X_{t_1}, {\mathcal{O}}_{X_{t_1}}(D_{t_1}+2cK_1) )
  \geq 
  h^0(X_{t}, {\mathcal{O}}_{X_{t}}(D_{t}+2cK_t) )=1,
  $$
 and   
  $$  h^0(X_{t_1}, {\mathcal{O}}_{X_{t_1}}(-D_{t_1}-2cK_1) )
  \geq 
  h^0(X_{t}, {\mathcal{O}}_{X_{t}}(-D_{t}-2cK_t) )=1,
  $$
  where $t\in B$, $D_t+2cK_t=cD_t' +2cK_t=c(D_t'+2K_t)=0$ (Lemma 2.2). So 
 ${\mathcal{O}}_{X_{t_1}}(D_{t_1}+2cK_1) $ must be trivial, i.e., 
 $$D_{t_1}+2cK_1=n(t_1)D_1'+2cK_1
 =-2n(t_1)K_1+2cK_1=0.$$ 
 Therefore  $n(t_1)=c$ for every $t_1\in V-B$.
 Hence $D_t=cD_t'$ for every $t\in V$.
 By changing
  coefficient locally, we may assume $D|_{X_t}=D_t'$, where $2D_t'+K_t=0$.

Since $2D_t+K_{X_t}=0$ for every $t\neq t_0$, similarly, considering  
the divisor 
$2D+K_X$, when  restricted to every fibre $X_t$,  
we have 
$$    h^0(X_0, {\mathcal{O}}_{X_0}(2D_0+K_0))
      \geq h^0(X_t, {\mathcal{O}}_{X_t}(2D_t+K_{X_t}))=1.
$$
On the other hand
$$    h^0(X_0, {\mathcal{O}}_{X_0}(-2D_0-K_0))
      \geq h^0(X_t, {\mathcal{O}}_{X_t}(-2D_t-K_{X_t}))=1.
$$
So we have 
$$h^0(X_0, {\mathcal{O}}_{X_0}(2D_0+K_0))
=h^0(X_0, {\mathcal{O}}_{X_0}(-2D_0-K_0))=1.$$
Again this  implies  that the sheaf
 ${\mathcal{O}}_{X_0}(2D_0+K_0))$ is trivial.
Hence  $2D_0+K_0=0$.
Since  $S_t$ has vanishing Hodge cohomology, and
 $X_0$ is isomorphic to $X_t$,
we have 
$$h^0(X_0, {\mathcal{O}}_{X_0}(2nD_0))
=h^0(X_0, {\mathcal{O}}_{X_0}(-nK_0))
$$
$$
=h^0(X_t, {\mathcal{O}}_{X_t}(-nK_t))
=h^0(X_t, {\mathcal{O}}_{X_t}(2nD_t'))=1.
$$
Therefore $S_0$ is not affine.
\begin{flushright}
 Q.E.D. 
\end{flushright}

\begin{remark} Let $U$ be covered by a set of small open discs $U_i$.
By the above argument, for each $i$, 
there is a constant  $c_i$ such that for $t\in U_i$,
$D_t=c_iD_t'$, where $D_i'$ is the irreducible boundary elliptic curve
on $X_t$.  Since $U$ is connected, all these $c_i's$ are equal. 
That is, there is constant $c$, such that for all $t\in U$,
$D|_{X_t}=D_t=cD_t'$.  Thus we have proved that 
by changing the coefficients of $D$,  the new boundary divisor $D'$
on $X$  satisfies $D'|_{X_t}=D_t'$.    
\end{remark}

\begin{lemma}
If there is an affine open set $U$ in $C$ such that 
for every $t\in U$, $t\neq t_0$, $S_t$ is of type (3) open surface, 
where $t_0$ is a fixed point of $U$, then $S_0$ must be of the same
type. 
\end{lemma}
$Proof$.               First, $S_0$ is not of type (2) open surface 
since $X_0$ is rational by 
   Iitaka's theorem [I4]. We only need to prove that 
$S_0$ is not affine as above lemma. It suffices to prove that
$h^0(X_0, {\mathcal{O}}_{X_0}(nD_0))<cn^2$ for all positive 
number $c$ (Lemma 2.1).  

 As in  Lemma 2.5, we may assume that $D_t$ 
contains no  exceptional curves of the first kind    
 for every
$t$ in $U$. In 
fact, if there is an exceptional curve $E_1$ of the first kind
in  $D_{t_1}$ for some point $t_1$ in $U$, then 
locally analytically,   $E_1$ sits in an 
irreducible nonsingular divisor $D_1$ of $X$, i.e., 
$D_1$ is a prime component of $D$. (We may assume that
$D$ is an effective divisor on $X$ with simple normal crossings \cite{Zh}.)
Now 
$f$ is proper  on $D_1$  
and $D_1$ is a manifold. So we can apply Kodaira's extension theorem 
locally on $D_1$ near $D_{t_1}$.
More precisely, in our case, we can compute it directly. 
 Since $D_1$ is smooth,
 for a small number $\epsilon>0$, in a neighborhood 
 $V=\{t\in C,  |t-t_1|< \epsilon\}$ of $t_1$,
 $D_1$    intersects every fibre $X_t$ 
 with  a prime  divisor on $X_t$. Since  
 $h^0({\mathcal{O}}_{X_t})=1$  and
 $h^1({\mathcal{O}}_{X_t})=
 h^2({\mathcal{O}}_{X_t})=0$ (Lemma 2.1, Lemma 2.2),
 by the Riemann-Roch   formula
and upper semi-continuity, we have 
$$ \chi ({\mathcal{O}}_{X_t}(nE_t))
=1+\frac{1}{2} n^2 E_t^2-
\frac{1}{2} n E_t\cdot K_{X_t}=1+\frac{1}{2} n^2 E_1^2-
\frac{1}{2} n E_1\cdot K_{X_1}, 
$$
where $E_1=D_1|_{X_t}$. 
So $E_t$ is again an $(-1)-$curve on $D_1$.  
This implies that all the extended $(-1)$ exceptional curves 
near $D_{t_1}$ sit in $D_1$ and do not meet $Y$. So
after contraction, $Y$ remains the same, that is, when contracting 
$(-1)-$curves, we only change the
boundary $D_t$  but all the  open surfaces $S_t$ over $t\in U$ 
remain unchanged.

If $D_t$ is the special divisor $D_t'$ as in Mohan Kumar's Theorem, i.e., 
its dual graph is either $\tilde{D}_8$ or $\tilde{E}_8$,
then 
 we have 
$D_t+K_{X_t}=0$ for every $t\in U$ and $t\neq t_0$ by Lemma 2.2.
 By the similar inequalities as in  the proof of Lemma 2.5,
we have  
$$h^0(X_0, {\mathcal{O}}_{X_0}(D_0+K_0))
=h^0(X_0, {\mathcal{O}}_{X_0}(-D_0-K_0))=1.$$
Hence  $D_0+K_0=0$.  Since $X_0$ is of type (3) projective 
surface 
 and $H^i(S_0, \Omega^j_{S_0})=0$,
we know  that  $S_0$ 
is not affine
and must be of   type (3) open surface. 
But we cannot guarantee 
that  the  dual graph of  $D_t$  is either $\tilde{D_8}$ or $\tilde{E_8}$.
 We only know $D_t$ has nine components
and every prime component is isomorphic to ${\Bbb{P}}^1$ 
with self-intersection $-2$  ([Ku]).
In   Lemma 2.5, we may  assume that the special divisor 
$D_t'$   on 
$X_t$ is the restriction 
of a global divisor $D$ on $X_t$ since $D_t'$  has only one component by
 Remark 2.6. 
Now the situation is more delicate.

  Let $D_t'$ be the special divisor of type (3) projective surface
  as  above, i.e., its dual graph is  either $\tilde{D_8}$ or $\tilde{E_8}$,
  $D_t'\cdot  D_t'=D_t'\cdot K_t=0$ and 
${\mathcal{O}}_{D_t'}(D_t')$ is nontorsion \cite{Ku}.
  For any nonnegative integer $n$, there is $m$ such that 
 $mD_t'-nD_t$ is effective. For example, we may choose 
 $m=an$ where $a$ is the maximum coefficient of $D_t$'s components. 
 So 
 $$0< h^0(X_t, {\mathcal{O}}_{X_t}(nD_t)) 
 \leq h^0(X_t, {\mathcal{O}}_{X_t}(mD_t'))=1.
 $$
 Therefore $h^0(X_t, {\mathcal{O}}_{X_t}(nD_t))=1.$ 
By Serre duality, $H^2(X_t, {\mathcal{O}}_{X_t}(nD_t))=0$
for all $n\gg 0$ since $K_{X_t}$  and $D_t$ have the same support
by Lemma 2.2. 
 Consider $h^1(X_t, {\mathcal{O}}_{X_t}(nD_t))$, there are three cases [Za].

 Case 1.  $h^1(X_t, {\mathcal{O}}_{X_t}(nD_t))$ is bounded, i.e., there is 
 a positive integer $k$ such that for all $n\geq 0$, we have
 $$
 h^1(X_t, {\mathcal{O}}_{X_t}(nD_t)) \leq k < \infty.
 $$
 By Zariski's theorem, page 611,  \cite{Za}, $D_t$ is arithmetically effective. 
 By Riemann-Roch formula and Lemma 2.1, Lemma 2.2,  we have
 $$
 h^1(X_t, {\mathcal{O}}_{X_t}(nD_t))=
 -\frac{1}{2}n^2D_t^2
 +\frac{1}{2}nD_t\cdot K_{X_t}.
 $$
 This equality gives us $D_t^2=D_t.K_t=0$ since 
 $h^1(X_t, {\mathcal{O}}_{X_t}(nD_t))$ is bounded for all $n$.
 Then for every prime component $E$ in $D_t$, we have 
 $E.D_t=0$ since $D_t$ is arithmetically effective. 
 By Lemma 1.7 [Ku], 
 we know $D_t=n(t)D_t'$, where the  positive integer $n(t)$ 
 depends on the point $t$ in $U$. So for every $n\geq 0$, 
by Lemma  2.1, we have
 $$
 h^1(X_t, {\mathcal{O}}_{X_t}(mD_t))=
 h^1(X_t, {\mathcal{O}}_{X_t}(mn(t)D_t'))=0.
 $$
 Now the Euler characteristic  of 
     ${\mathcal{O}}_{X_0}(nD_0)$ is
$$
\chi({\mathcal{O}}_{X_0}(nD_0))
=1- \frac{1}{2}n^2D_0^2 +
\frac{1}{2}nD_0.K_{X_0}
=1.
$$   
Thus  again $D_0^2=D_0\cdot  K_{X_0}=0$. By the same argument 
as in the above 
lemma and remark,
there is a positive integer $c$ such that for every $t\in U$ and $t\neq t_0$,
$D_t=cD_t'$. Considering  the divisor $D+cK_X$ on $X$, when restricted 
to $X_0$, we have 
$$h^0(X_0, {\mathcal{O}}_{X_0}(D_0+cK_{X_0}))\geq 1, \quad
h^0(X_0, {\mathcal{O}}_{X_0}(-D_0-cK_{X_0}))\geq 1.$$
 So
$D_0=-cK_{X_0}$. 
Let $D_0=P+N$ be the Zariski decomposition of $D_0$,
then $P$ is nef, $N$ is  negative  definite (both are effective)
and every 
component of $N$ does not intersect $P$.  Let $E$ be a prime component
of $P$. Locally analytically, $E$ is contained in a prime 
divisor $G$ of $X$. Let $G|_{X_t}=E_t$. Applying  upper 
semi-continuity and Riemann-Roch formula
to ${\mathcal{O}}_{X_0}(nE)$ and 
${\mathcal{O}}_{X_t}(nE_t)$, we have 
$E\cdot K_0=E_t\cdot K_t=0$, \cite{Ku}, Lemma 3.1.  Thus 
 $E \cdot D_0= E \cdot (-cK_0)=0$. 
If $E\cdot P>0$, 
then $E\cdot P=  E\cdot (D_0-N)=-E\cdot N>0$.
Therefore $E\cdot N<0$. This means that $E$ must be a component of $N$
which  is a contradiction since no component of $N$ intersects $P$.
So $P^2=0$. By Corollary 14.18,  \cite{Ba}, 
$\kappa(D_0, X_0 )\leq 1$. By Lemma 2.1, 
$S_0$ is not affine.

Case 2. If  $h^1(X_t, {\mathcal{O}}_{X_t}(nD_t))$ is as large as $cn$
for some positive number $c$, then by Zariski's theorem [Za], page 611, 
$D_t$ is arithmetically effective and 
the intersection form of $D_t$ is negative definite. This 
contradicts  Lemma 1.6 [Ku]. So this case cannot happen.

Case 3. If  $h^1(X_t, {\mathcal{O}}_{X_t}(nD_t))$ is as large as $kn^2$
for some positive number $k$, then by Riemann-Roch formula,
we know  $D_t^2<0$. Let $D_t=A+B$ be the Zariski decomposition such that 
$A$ is arithmetically effective, $B\geq 0$ is negative definite
and every prime component of $B$ does not meet $A$. There is a positive 
integer $n_0$ such that $n_0A$ and $n_0B$ are integral. Without loss 
of generality, we may assume $A$ and $B$ are integral. Since there is  a
positive integer $l$
such that $lD_t'-D_t$ is effective, we have exact sequence
$$
0\longrightarrow {\mathcal{O}}(nD_t)
\longrightarrow  {\mathcal{O}}(nlD_t')
\longrightarrow  Q
\longrightarrow 0,
$$
where $Q$ is the cokernel. Hence we still have
 $h^0(X_t, {\mathcal{O}}_{X_t}(nD_t))=1$  even though 
$D_t$ is different from $D_t'$. Therefore $\kappa(D_t, X_t)=0$
by Iitaka's Theorem.
This implies $A^2=0$ \cite{Za} or 
\cite{Ba}, Corollary 14.18.  
Since $A$ is arithmetically effective
and supp$D_t=$supp$A\cup$supp$B$, for every prime component $E$
of $D_t$, $E.A=0$. By Corollary 1.7 [Ku], there is positive integer 
$m_0$ such that $A=m_0D_t'$. So $D_t^2=B^2$ and $D_t-D_t'\geq 0$.
Let $D_{0,i}$ be a prime component of $D_0=D|_{X_0}$.
 Choose a small neighborhood
$V$ of $t_0$ such that locally analytically in $V$, $D_{0,i}$ lies 
in a unique
prime divisor $D_i$ of $f^{-1}(V)$. $D_i$ cuts every fibre $X_t$, $t\in V$
with an irreducible  $(-2)$-curve. So over $V$, there is one to one 
correspondence between the prime divisor of $D_t$ and the prime divisor 
of  $f^{-1}(V)$. We may rearrange the coefficients of $D_i$ 
locally as in the proof of the above lemma such that 
$D_t=cD_t'$, $t_0\neq t\in V$. So 
$h^1(X_t, {\mathcal{O}}_{X_t}(nD_t))=0$ for all 
$t\in V$,  $t\neq t_0$. 
Then we reduce Case 3 to Case 1. This proves that $S_0$ is not affine.
\begin{flushright}
 Q.E.D. 
\end{flushright}

 \begin{remark}
If $D_t$ is not the special divisor as in Mohan Kumar's Theorem, i.e., 
$D_t$ has different coefficients from $D_t'$ but they have the same 
support, then
we still have $D_t\cdot K_t=0$ by Lemma 3.1 \cite{Ku}. 
Since $h^0(X_t, {\mathcal{O}}_{X_t}(nB))=1$,
by Riemann-Roch formula,
$h^1(X_t, {\mathcal{O}}_{X_t}(nD_t))=-\frac{1}{2}n^2D_t^2
=-\frac{1}{2}n^2B^2
=h^1(X_t, {\mathcal{O}}_{X_t}(nB))\sim cn^2.$
Thus $B^2<0$. Since  $B$ is definite negative,
by  Lemma 1.6 [Ku], we know that
the support of $B$ is strictly smaller than 
the support of $D_t$.   \\
\end{remark}

 \begin{remark} Let $X_t$ be of type (3) projective surface  and  
$D_t'$ the special 
divisor  as above. Let $E$ be any prime component of $D_t'$.  Then $E^2=-2$
[Ku]. Since the canonical divisor $K_{t}=-D_t'$, by Riemann-Roch, 
$h^1(X_t, {\mathcal{O}}_{X_t}(nD_t'+E))=n^2$. 
Combining with the above argument, 
for any divisor $D_t$ with the same support as $D_t'$, we know either 
$h^1(X_t, {\mathcal{O}}_{X_t}(nD_t))=0$ or 
$h^1(X_t, {\mathcal{O}}_{X_t}(nD_t))\sim cn^2$ for some positive integer $c$. 
\end{remark}

\begin{lemma}If $S_0$ is affine then  
there is an affine open set $U$ in $C$ such that 
for every $t\in U$, $S_t$ is affine.
\end{lemma}
$Proof.$  This is a direct consequence of the above lemmas
since $S_0$ only can be one of the three surfaces. 
\begin{flushright}
 Q.E.D. 
\end{flushright}

\begin{lemma} If there is an affine open set $U$ in $C$ such that 
for every $t\neq t_0$, $S_t$ is affine, then $S_0$ is affine.
\end{lemma}
$Proof.$ This is an immediate conclusion of Mohan Kumar's 
classification and upper
semi-continuity theorem. 
\begin{flushright}
 Q.E.D. 
\end{flushright}

The first half of Theorem 1.1 follows from the above lemmas. 
The second half is a direct
consequence 
of Theorem 5.11 and Theorem 6.12
in \cite{U}. In fact, if one smooth fibre $X_0$ is not affine, 
then all smooth fibres are not affine by the above lemmas. 
Since $X_0$ is a ruled surface, $\kappa(X_0)=-\infty.$
So 
$$\kappa(X)\leq \kappa(X_0)+1=-\infty.$$ 
By Lemma 2.1,   
$$  0< \kappa(D, X) \leq \kappa(D_t, X_t)+1=1.
$$
This completes the proof of Theorem 1.1.

%\tableofcontents
%%
%%
%%
%%
%%
%%
%%
\section{Proof of Theorem 1.2}

\begin{lemma}[\bf{Goodman, Hartshorne}]
{\it  Let $V$ be a scheme and  $D$ be an 
effective Cartier divisor on  $V$. Let $U=V-{\Supp}D$ and $F$ 
be any coherent sheaf on $V$, 
then for every $i\geq 0,$} 
\end{lemma}
$$\lim_{{\stackrel{\to}{n}}}
H^i(V, F\otimes {\mathcal{O}}(nD)) \cong  H^i(U,  F|_U).
$$

This  lemma enables  us  to transfer  the cohomology 
information  from $Y$  to its completion $X$.

$Proof$  $of$   $Theorem$  $1.2$. The idea is to prove  
for any coherent sheaf 
$F$ on $Y'$, $H^i(Y', F_{Y'})=0$ for all $i>0$. 
Since the dimension of $Y'$ is 3, we only need to 
consider $i=1,2,3.$
We use the technique
 in \cite{Zh} with some modification. We present all details 
for completeness.

Notice that  $Y'\subset Y$.
Let ${{F}}_{Y'}$ be any 
coherent sheaf on $Y'$, 
then it can be extended to a coherent sheaf ${{F}}_X$
on $X$ and ${{F}}_Y|_{S_t}$, ${{F}}_X|_{X_t}$
are coherent, \cite{H1}  page 115, page 126. We will not distinguish 
them and just write ${{F}}$. 
Since general fibre $X_t$ over $t\in C$ is smooth and irreducible
\cite{Zh} and for any  $F$, there is an open set $U$ in $C$ such that
$R^if_*F$ is locally free on $U$, we may assume that $R^if_*F$ is locally
free on $C$ and every fibre over $C$ is smooth and irreducible.

(1) Proof of  $H^3(Y, F)=0$.

Since $S_t$ is affine for every
$t$ in $C$, we have $H^i(S_t, {{F}}|_{S_t})=0$ for every 
$i>0$. Let 
$F_n=F\otimes {\mathcal{O}}_X(nD)$ and 
$F_{n,t}=F\otimes {\mathcal{O}}_X(nD)|_{X_t}$.
 By Goodman and Hartshorne's Lemma, 
$$\lim_{{\stackrel{\to}{n}}}
H^i(X_t, F_{n,t})=0
$$
for all $i>0$ and $t\in C$. 
Since each fibre has dimension 2, we have $H^3(X_t, F_{n,t})=0$
for all $n\geq 0$ and $t\in C$.   By upper semi-continuity,
$R^3f_*F_n=0$ for all $n$. Again  by Goodman and Hartshorne's Lemma,
$$ H^3(Y, F)=\lim_{{\stackrel{\to}{n}}}
H^3(f^{-1}(C), F_n)=
\lim_{{\stackrel{\to}{n}}}
R^3f_*F_n(C)=0.
$$

(2) Proof of  $H^2(Y, F)=0$.

It suffices to prove  the claim for locally free sheaves. 
In fact, suppose  $H^2(Y, L)=0$ for any locally free sheaf $L$ on $X$.
For any coherent sheaf $F$ on $X$, there is a locally free sheaf 
$L$ on $X$ such that we have the surjective map 
$L\longrightarrow F$. Let $K$ be the kernel, then we have short 
exact sequence  on $Y$
$$ 0\longrightarrow 
K
\longrightarrow  L
\longrightarrow  F
\longrightarrow 
0.
$$
By step 1, we know $H^3(Y, K )=0$ since $K$ is also coherent \cite{H1}. So 
$H^2(Y, L)=0$ implies $H^2(Y, F)=0$.

So we may assume that $F$ is a locally free sheaf on $X$.

Let $t\in C$.
 From the exact sequence
$$0\longrightarrow 
{\mathcal{O}}_X(nD)
\longrightarrow
{\mathcal{O}}_X((n+1)D) 
\longrightarrow
{\mathcal{O}}_D((n+1)D) 
\longrightarrow
0, 
$$
tensoring with $F$ then with ${\mathcal{O}}_{X_t}$, 
  we have 
$$0\longrightarrow 
F_{n,t}
\longrightarrow
F_{n+1, t}
\longrightarrow
F_{n+1, t}|_{D_t}
\longrightarrow
0. 
$$
Since $D_t$ is a curve, 
$H^2(X_t, F_{n+1, t}|_{D_t})=0$ for all 
$ n\geq 0$ and $t\in C$. 
So the map $H^2(X_t, F_{n, t} )
\rightarrow H^2(X_t, F_{n+1, t} )$
is surjective. Since $S_t=X_t-D_t$ is affine, by Goodman and Hartshorne's 
Lemma,
$$\lim_{{\stackrel{\to}{n}}}
H^2(X_t, F_{n,t})
=H^2(S_t, F )=0.
$$  
So for any $t\in C$, there is a positive integer
$n(t)$,  depending on $t$ such that for every 
$n\geq n(t)$, 
$H^2(X_t, F_{n,t} )=0$.

Given any $n$, there is  an  affine open set $U_n$
of $C$ such that $R^2f_*F_n$  is locally free 
on $U_n$. By the same argument as in the  next paragraph,
the intersection of these infinitely many open sets 
is   not empty. 
 Now fix some $t_0$ in $\cap U_n$ such that 
$H^2(X_{t_0}, F_{n, t_0})=0$ for every $n\geq n(t_0)$ 
and there is an open neighborhood
 $U_0$ of $t_0$ in $\bar{C}$ such that $R^2f_*F_{n(t_0)}$  is locally free  
on   $U_0$. Then $H^2(X_t, F_{n(t_0), t})=0$ for every $t$ in  $U_0$. So   
$H^2(X_t, F_{n, t})=0$ for every $t$ in  $U_0$ and every $n\geq n(t_0)$. 
Let $C-U_0=\{t_1, t_2, ..., t_k\}$, choose 
$n_0=\max (n(t_0), n(t_1), ..., n(t_k))$, 
then  
$H^2(X_t, F_{n, t})=0$
for every  $t\in C$ and every $n\geq n_0$.
By  upper semi-continuity theorem,  
$(R^2f_*F_n)_t/{\mathcal{P}}(R^2f_*F_n)_t=0$
for all points $t$ in $C$. By  Nakayama's lemma,   $R^2f_*F_n|_C=0$. 
By Goodman and Hartshorne's Lemma
$$ H^2(Y, F)=
\lim_{{\stackrel{\to}{n}}}
H^2(f^{-1}(C), F_n)=
\lim_{{\stackrel{\to}{n}}}
R^2f_*F_n(C)=0.
$$

(3) Proof of  $H^1(Y', F)=0$,  where $Y'$ is an open subset of 
$Y$ obtained by removing finitely many fibres from $Y$.

Let $F_n$ be as above. 
For any fixed $n$, 
there is an open set $U_n$ in  
 $\bar{C}$, such that $R^1f_*F_n$ is locally free on  $U_n$. Let 
$U_n=\bar{C}\backslash A_n$, where $A_n$ is closed in  $\bar{C}$, i.e., 
it consists only finitely many points of  $\bar{C}$.  Since any complete 
metric space is a Baire space, 
\cite{Bo2},    Chapter 9, 
$B=\bar{C}\backslash \cup A_n=\cap U_n$ is a dense 
subset of $\bar{C}$ in complex topology. Hence for every point $t$ in $B$,
 all stalks    
$(R^1f_*F_n)_t$ are locally free. Write $B$ as a union of connected subsets 
$B_m$, $B=\cup B_m$, then there is one  $B_m$, such that $B_m$ is dense in 
$\bar{C}$ and connected in complex topology. So we may assume that $B$ 
is connected.  
Again by upper-semicontinuity theorem, for every point $t$ in $C$
 and every $n\geq n_0$, since  $R^2f_*F_n|_C=0$, we have  
$$(R^1f_*F_n)_t{\otimes}
{\Bbb{C}} \cong H^1(X_t, F_{n, t}).
$$
For any  $m$, $h^1(X_t, F_{m, t})$ is constant on $B$ since 
$R^1f_*F_m$ is locally free at every point $t$ in $B$
and $B$ is connected. So for the above $n$ and for all points  
$t$ in  $B$, there is $l$ such that  the map 
$$
H^1(X_t, F_{n, t})\longrightarrow H^1(X_t, F_{n+l, t})
$$
is zero.  Moreover, for every point $t$ in $C$ and sufficiently large 
$n$, we have the following 
commutative diagram 

\[
  \begin{array}{ccc}
    R^1f_*F_n{\otimes}
{\Bbb{C}}(t)                            &
     {\stackrel{\approx}{\longrightarrow}} &
     H^1(X_t, F_{n, t})                               \\
    \Big\downarrow\vcenter{%
        \rlap{$\scriptstyle{\alpha}$}}              &  &
    \Big\downarrow\vcenter{%
       \rlap{$\scriptstyle{\beta}$}}      \\
 R^1f_*F_{n+l}{\otimes}
{\Bbb{C}}(t)             &  {\stackrel{\approx}{\longrightarrow}}  &
 H^1(X_t, F_{n+l, t}). 
\end{array}
\]
The map  $\beta$ is zero for every $t\in B$, so  the map
$$\alpha :\quad
(R^1f_*F_n)_t/{\mathcal{P}}(R^1f_*F_n)_t
\longrightarrow 
(R^1f_*F_{n+l})_t/{\mathcal{P}}(R^1f_*F_{n+1})_t
$$
is zero for all  points $t$ in $B$.
By the local freeness, this says  for every point $t$ in $B$, 
the stalks satisfy   
$$\lim_{{\stackrel{\to}{n}}}
(R^1f_*F_n)_t=0.  
$$ 
To see this, 
fix a point $t_0$ in $B$, for any sufficiently large $n$
and for the above $l$, choose an affine open set $V$ containing 
$t_0$ such that both $R^1f_*F_n$ and 
 $R^1f_*F_{n+l}$ are locally free on $V$. So there are two 
 positive integers $m_1$ and $m_2$ such that 
$R^1f_*F_n(V)= {\mathcal{O}}(V)^{m_1}$ 
and 
$R^1f_*F_{n+l}(V)= {\mathcal{O}}(V)^{m_2}$. Now for infinitely many 
maximal ideals ${\mathcal{P}}$,  we have commutative 
diagram

\[
  \begin{array}{ccc}
    {\mathcal{O}}(V)^{m_1}                            &
     {\stackrel{\psi}{\longrightarrow}} &
   {\mathcal{O}}(V)^{m_2}                                \\
    \Big\downarrow\vcenter{%
        \rlap{$\scriptstyle{\pi_1}$}}              &  &
    \Big\downarrow\vcenter{%
       \rlap{$\scriptstyle{\pi_2}$}}      \\
  {\mathcal{O}}(V)^{m_1} /{\mathcal{P}}{\mathcal{O}}(V)^{m_1}              &  
  {\stackrel{\phi}{\longrightarrow}}  &
  {\mathcal{O}}(V)^{m_2}/{\mathcal{P}}{\mathcal{O}}(V)^{m_2}.   
\end{array}
\]

Since  $\psi({\mathcal{O}}(V)^{m_1})
\subset \cap {\mathcal{P}}{\mathcal{O}}(V)^{m_2}=0$, 
where $\mathcal{P}$ runs over infinitely many maximal ideals of 
${\mathcal{O}}(V)$, we have $\psi({\mathcal{O}}(V)^{m_1})=0$. This proves 
$$\lim_{{\stackrel{\to}{n}}}
(R^1f_*F_n)_t =0.  
$$  

Since the direct limit of $R^1f_*F_n$  
 is quasi-coherent, 
its support is locally closed.  Now
 $B$ is dense and connected in complex topology, 
there exists an affine open set $U$ 
in $\bar{C}$
such that on $U$, the direct limit  
$$\lim_{{\stackrel{\to}{n}}}
R^1f_*F_n|_U =0.  
$$
Let $Y'=f^{-1}(U)-D$, 
by Goodman and Hartshorne's Lemma,
we have
$$H^1(Y', F)=
\lim_{{\stackrel{\to}{n}}}
H^1(f^{-1}(U), F_n)
=\lim_{{\stackrel{\to}{n}}}
R^1f_*F_n(U)=0. 
$$
This finishes the proof of Theorem 1.2.

\begin{remark} In our  above proof of step 3,
 we encounter the following two  
questions if we do not know the local freeness  of $R^1f_*F_n$. 

(1) Let  $U$ be a smooth affine curve, then ${\mathcal{O}}(U)=A$
is a Dedekind domain. Let $N$ be a finitely generated module over 
$A$, then under what condition, $\cap ({\mathcal{P}}N)=0$?
Where ${\mathcal{P}}$ runs over all maximal ideals of $A$. 
A sufficient condition is that  $N$ is projective module. 
But this condition is too strong. Our $N$ is defined by cohomology.
 It is hard to see it is projective or not. Definitely 
finitely generated module  is not sufficient. For example,
let $U={\Bbb{A}}^1$, $A={\mathcal{O}}(U)={\Bbb{C}}[x]$,
$N={\Bbb{C}}[x]/(x^2)$, then  $\cap ({\mathcal{P}}N)\neq 0.$  

(2) Let  $A$ be a Dedekind domain  and $\mathcal{P}$ as above, 
let $(M_n, f_n)$
be a direct system of finitely generated $A$-modules.
If 
$$ 
\lim_{{\stackrel{\to}{n}}} (M_n/{\mathcal{P}} M_n)=0,
$$
under what conditions, can  we  say that 
$$\lim_{{\stackrel{\to}{n}}}M_n=0?
$$
Again all $M_n$ being finitely generated is not sufficient. 
For example, let $A={\Bbb{C}}[[t]]$, the ring of formal 
power series, let $M_n=t^{-n}A$, then 
$$\lim_{{\stackrel{\to}{n}}}M_n={\Bbb{K}}\neq 0, 
$$
where ${\Bbb{K}}= {\Bbb{C}}((t))$ but
$$\lim_{{\stackrel{\to}{n}}}M_n/{\mathcal{P}} M_n=0.
$$
  \\

\end{remark}

%\tableofcontents
%%
%%
%%
%%
%%
%%
%%
\section{Proof of Theorem 1.3}

\begin{lemma} $R^if_*{\mathcal{O}}_X(nD)$ is locally free for all $i\geq 0$
and $n\gg 0$.
\end{lemma}
$Proof$.  Since each fibre has dimension 2, by upper semi-continuity
theorem, $R^if_*{\mathcal{O}}_X(nD)=0$ for all $i>2$ and
$n\geq 0$. By Lemma 2.1, since each fibre  $X_t$
is of type (2), we have 
$$h^0(X_t, {\mathcal{O}}_{X_t}(nD_t))
=h^1(X_t, {\mathcal{O}}_{X_t}(nD_t))=1 \quad
\mbox{and} \quad 
h^2(X_t, {\mathcal{O}}_{X_t}(nD_t))=0$$ for all $t\in C$ and $n\gg 0$.  
\begin{flushright}
 Q.E.D. 
\end{flushright}

\begin{lemma} $R^if_*\Omega^3_X(nD)=R^if_*{\mathcal{O}}_X(K_X+nD)$ 
is locally free for all $i\geq 0$ and $n\gg 0$.
\end{lemma}
$Proof$. Since $X_t$ is smooth, we have $K_X+D|_{X_t}=K_{X_t}=K_t$. 
So $\Omega^3_X(nD)|_{X_t}=
{\mathcal{O}}_X(K_X+nD)|_{X_t}=
{\mathcal{O}}_{X_t}(K_t+(n-1)D_t),$
where $D_t=D|_{X_t}.$
By Lemma 2.1, 
$$h^0(X_t, \Omega^3_X(nD)|_{X_t})=
h^0(X_t, {\mathcal{O}}_{X_t}(K_t+(n-1)D_t)= 1,$$
 $$h^1(X_t, \Omega^3_X(nD)|_{X_t})=
h^1(X_t, {\mathcal{O}}_{X_t}(K_t+(n-1)D_t)=1,$$
and  
$$h^2(X_t, \Omega^3_X(nD)|_{X_t})=
h^2(X_t, {\mathcal{O}}_{X_t}(K_t+(n-1)D_t)=0$$
for all  $n\gg 0$.  This proves  the local freeness. 
\begin{flushright}
 Q.E.D. 
\end{flushright}

\begin{lemma} $R^if_*\Omega^1_X(nD)$ 
is locally free for all $i\geq 0$ and $n\gg 0$.
\end{lemma}
$Proof$. From the exact sequences  ([H1], II, Theorem 8.17 and 
 [GrH], page 157)
 $$0\longrightarrow 
 {\mathcal{O}}_{X_{t}}
 \longrightarrow \Omega^1_X|_{X_t}
 \longrightarrow \Omega^1_{X_t}
 \longrightarrow  0, 
 $$
 tensoring with ${\mathcal{O}}_X(nD)$, we have
  $$0\longrightarrow 
 {\mathcal{O}}_{X_{t}}(nD_t)
 \longrightarrow \Omega^1_X(nD)|_{X_t}
 \longrightarrow \Omega^1_{X_t}(nD_t)
 \longrightarrow  0.
 $$
We will prove that for any two points $t, t'\in C$ and all $n\gg 0$,
we have  
$$h^i(X_t, \Omega^1_X(nD)|_{X_t})=h^i(X_{t'}, \Omega^1_X(nD)|_{X_{t'}}).$$
Then by upper semi-continuity theorem, we are done. 
By the above short exact sequences,  for both fibres $X_t$ and $X_{t'}$,
we have the commutative diagram
\[
  \begin{array}{ccccccccc}
0\longrightarrow H^0({\mathcal{O}}_{X_t}(nD_t))&
{\stackrel{\alpha_1}{\longrightarrow}}  H^0(\Omega^1_X(nD)|_{X_t})&
{\stackrel{\alpha_2}{\longrightarrow}}  H^0(\Omega^1_{X_t}(nD_t))&
{\stackrel{\alpha_3}{\longrightarrow}}   H^1({\mathcal{O}}_{X_t}(nD_t))\\
\quad \quad   \parallel &
\Big\downarrow\vcenter{%
        \rlap{$\scriptstyle{\phi}$}} &
\quad \quad   \parallel &
\quad \quad   \parallel & \\
 0\longrightarrow H^0({\mathcal{O}}_{X_{t'}}(nD_{t'}))&
{\stackrel{\beta_1}{\longrightarrow}}  H^0(\Omega^1_X(nD)|_{X_{t'}})&
{\stackrel{\beta_2}{\longrightarrow}} H^0(\Omega^1_{X_{t'}}(nD_{t'}))&
{\stackrel{\beta_3}{\longrightarrow}}  H^1({\mathcal{O}}_{X_{t'}}(nD_{t'})),        
\end{array}
\]
where the natural map $\phi$ is defined as follows. If 
$\xi\in H^0(\Omega^1_X(nD)|_{X_t})$ is contained in the image of 
$H^0({\mathcal{O}}_{X_t}(nD_t))=\Bbb{C}$, then there is  a number
$a\in \Bbb{C}$ such that  $\xi=\alpha_1(a)$.  Thus we define
  $\phi(a)=\beta_1(a).$ 
If $\xi$ is not contained in the image of $\alpha_1$, then
$\alpha_2(\xi)\in H^0(\Omega^1_{X_t}(nD_t))$  and
$\alpha_3\circ \alpha_2(\xi)=0$. So there is 
$\eta\in H^0(\Omega^1_X(nD)|_{X_{t'}})$ 
such that  $\beta_3\circ \beta_2(\eta)=0.$
Define  $\phi (\xi)=\eta$, by 5-Lemma \cite{La},  we have  
$$H^0(\Omega^1_X(nD)|_{X_t})=
H^0(\Omega^1_X(nD)|_{X_{t'}}).
$$   
Similarly, we have 
$$H^i(\Omega^1_X(nD)|_{X_t})=H^i(\Omega^1_X(nD)|_{X_{t'}})$$ 
for $i>0$. 
\begin{flushright}
 Q.E.D. 
\end{flushright}

\begin{lemma} $R^if_*\Omega^2_X(nD)$ 
is locally free for all $i\geq 0$ and $n\gg 0$.
\end{lemma}
$Proof.$ Notice that we have the short exact sequence
(\cite{H1}, II, Theorem 8.17 and 
 \cite{GrH}, page 157)
 $$ 
 0\longrightarrow \Omega^1_{X_t}
 \longrightarrow \Omega^2_X|_{X_t}
 \longrightarrow \Omega^2_{X_t}
 \longrightarrow 0.
 $$
 Tensoring with ${\mathcal{O}}_X(nD)$, we have
$$ 
 0\longrightarrow \Omega^1_{X_t}(nD_t)
 \longrightarrow \Omega^2_X(nD)|_{X_t}
 \longrightarrow \Omega^2_{X_t}(nD_t)
 \longrightarrow 0.
 $$
By Lemma 2.1, for every $t\in C$,  we know
$$ h^0(X_t,\Omega^2_{X_t}(nD_t))=
h^1(X_t, {\mathcal{O}}_{X_t}(K_t+nD_t))=1\quad \mbox{and}\quad
h^2(X_t, {\mathcal{O}}_{X_t}(K_t+nD_t))=0.
$$

By the same argument as the above lemma, for any two points 
 $t, t'\in C$ and all $n\gg 0$, writing the long exact sequences
for $t$ and $t'$, 
we have  
$$h^i(X_t, \Omega^2_X(nD)|_{X_t})=h^i(X_{t'}, \Omega^2_X(nD)|_{X_{t'}}).$$
\begin{flushright}
 Q.E.D. 
\end{flushright}

\begin{lemma} For every $t\in C$, 
$H^i(S_t, \Omega^j_Y|_{S_t})=0$  
for all $j\geq 0$ and $i>0$, where  $S_t=X_t-D_t$.
\end{lemma}
$Proof.$ Since $S_t$ is a surface, we only need to consider 
$i=1,2.$ 

The   claim is obvious for ${\mathcal{O}}_Y$.  Since $S_t$ is smooth,
we  have  the  exact sequence
$$0\longrightarrow 
 {\mathcal{\phi}}_t/{\mathcal{\phi}}_t^2={\mathcal{O}}_{S_t}
 \longrightarrow 
 \Omega^1_Y|_{S_t}
 \longrightarrow
\Omega^1_{S_t} 
\longrightarrow
 0,
$$
where ${\mathcal{\phi}}_t$ is the defining sheaf of $S_t$. 
Therefore the claim holds for $\Omega^1_Y|_{S_t}$. 
Since the normal sheaf 
$${\mathcal{N}}_{S_t/Y}=
{\mathcal{H}}om({\mathcal{\phi}}_t/{\mathcal{\phi}}_t^2,  
{\mathcal{O}}_{S_t})
={\mathcal{O}}_{S_t},
$$
we have  
$${\mathcal{\omega}}_{S_t}
\cong {\mathcal{\omega}}_Y
\otimes 
{\mathcal{N}}_{S_t/Y}
\cong 
{\mathcal{\omega}}_Y
\otimes 
{\mathcal{O}}_{S_t}
= {\mathcal{\omega}}_Y|_{S_t}. 
$$
Hence the claim holds for $\Omega^3_Y|_{S_t}$.
From the exact sequence 
$$
0\longrightarrow 
\Omega^1_{S_t}
\longrightarrow 
\Omega^2_Y|_{S_t} 
\longrightarrow 
\Omega^2_{S_t} 
\longrightarrow 
0,
$$
we get the claim for  $\Omega^2_Y|_{S_t}$. 
\begin{flushright}
 Q.E.D. 
\end{flushright}

\begin{lemma} For all $j\geq 0$, $H^2(Y, \Omega_Y^j)=0$.
\end{lemma}
$Proof.$  The sheaves
$\Omega^j_X\otimes {\mathcal{O}}_X(nD)$
 are locally free by the above lemmas
for all $n\gg 0$.
 By upper semi-continuity theorem,
  there is an integer $n_0$ such that for all 
$n\geq n_0$, $R^2f_*\Omega^j_X\otimes {\mathcal{O}}_X(nD)|_C=0 $. 
By Goodman and Hartshorne's Lemma \cite{GH},
we have 
$$H^2(Y, \Omega_Y^j)=
\lim_{{\stackrel{\to}{n}}}
H^2(f^{-1}(C), \Omega^j_X\otimes {\mathcal{O}}_X(nD))=
\lim_{{\stackrel{\to}{n}}}
R^2f_*\Omega^j_X\otimes {\mathcal{O}}_X(nD)(C)=0.
$$
\begin{flushright}
 Q.E.D. 
\end{flushright}

\begin{lemma} For all $j\geq 0$, $H^1(Y, \Omega_Y^j)=0$.
\end{lemma}
$Proof.$   By the local freeness lemmas and Goodman and 
Hartshorne's Lemma, we have
 $$H^1(Y, \Omega_Y^j)=
\lim_{{\stackrel{\to}{n}}}
H^1(f^{-1}(C), \Omega^j_X\otimes {\mathcal{O}}_X(nD))=
\lim_{{\stackrel{\to}{n}}}
R^1f_*\Omega^j_X\otimes {\mathcal{O}}_X(nD)(C)=0.
$$
\begin{flushright}
 Q.E.D. 
\end{flushright}

%\tableofcontents
%%
%%
%%
%%
%%
%%
%%
\section{Proof of Theorem 1.6}

 We will prove Theorem 1.6 by construct an example. 
Let $C_t$ be a smooth projective elliptic curve 
defined by $y^2=x(x-1)(x-t)$, $t\neq 0, 1$.  Let $Z$ be the elliptic surface
defined by the same equation, then we have surjective morphism from $Z$ to 
$C={\Bbb{C}}-\{0, 1\}$ such that for every $t\in C$, the fibre $f^{-1}(t)=C_t$.   

\begin{lemma}
    There is a rank 2 vector bundle $E$ on $Z$ such that when 
   restricted to $C_t$, $E|_{C_t}=E_t$ is the unique nonsplit 
   extension of 
    ${\mathcal{O}}_{C_t}$ by ${\mathcal{O}}_{C_t}$, where $f$
     is the morphism from $Z$ to $C$.  
\end{lemma}    
    $Proof$.  Since $f:Z\longrightarrow C$ is an  elliptic fibration,
    for every $t$, we have   
$$h^1(f^{-1}(t), {\mathcal{O}}_{f^{-1}(t)})
=h^1({\mathcal{O}}_{C_t})=1.
$$
So
$$R^1f_*{\mathcal{O}}_Z{\otimes}
{\Bbb{C}}(t) \cong H^1(C_t, {\mathcal{O}}_{C_t})
\cong
{\Bbb{C}}.
$$
It gives us 
$$(R^1f_*{\mathcal{O}}_Z)_t/{\mathcal{P}}_t
(R^1f_*{\mathcal{O}}_Z)_t
\cong
{\Bbb{C}}.
$$
By Nakayama's lemma, $R^1f_*{\mathcal{O}}_Z$
is a line bundle on $C$. Since ${\Bbb{C}}[x, 1/x, 1/(x-1)]$
is principle ideal domain, the Picard group of $C$ is trivial,
i.e., any line bundle on $C$ is trivial. Therefore 
$R^1f_*{\mathcal{O}}_Z
\cong {\mathcal{O}}_C$
and 
$$H^1(Z, {\mathcal{O}}_Z)
=R^1f_*{\mathcal{O}}_Z(C)
={\mathcal{O}}_C(C)
={\Bbb{C}}[x, \frac{1}{x}, \frac{1}{x-1}].
$$
Given any exact sequence of vector bundles
$$  0\longrightarrow 
{\mathcal{O}}_Z
\longrightarrow 
E
\longrightarrow 
{\mathcal{O}}_Z
\longrightarrow 
0,
$$
let  $\xi$  be the image of  unit of  
$H^0(Z, {\mathcal{O}}_Z)$
in $H^1(Z, {\mathcal{O}}_Z)$, we get an element 
of $H^1(Z, {\mathcal{O}}_Z)$.  Conversely,
given any element   $\xi$ in $H^1(Z, {\mathcal{O}}_Z)$,
we can get an exact sequence as above. The procedure is the following.
Take any  (degree) large  ample line bundle $L$  on the elliptic 
surface $Z$,
for any positive integer $n$, 
we have   an exact sequence
$$  0\longrightarrow 
{\mathcal{O}}_Z
{\stackrel{\alpha }{\longrightarrow}}
L^{\oplus n}
{\stackrel{\beta }{\longrightarrow}}
G
\longrightarrow 
0,
$$ 
where  $G$ is  the quotient  which is a vector bundle.  We may assume  
$H^1(Z, L)=0$ by  raising the  degree  of $L$  since $L$ is  ample. 
So we have surjective map 
$H^0(Z, G)\twoheadrightarrow H^1(Z, {\mathcal{O}}_Z )$.
Hence  $\xi$  can be lifted to an element  $\eta$ in 
$H^0(Z, G)$.    
This element  $\eta$  defines  a map  
from ${\mathcal{O}}_Z$ to $G$,  
$\eta: {\mathcal{O}}_Z \rightarrow  G$,
sending 1 to $\eta$.  
Let  $E=\beta ^{-1}(\eta ({{\mathcal{O}}_Z}))$,
then we have exact sequence 
$$  0\longrightarrow 
{\mathcal{O}}_Z
{\stackrel{\alpha }{\longrightarrow}}
E
{\stackrel{\beta }{\longrightarrow}}
\eta ({{\mathcal{O}}_Z})={\mathcal{O}}_Z
\longrightarrow 
0.
$$
So there is one-to-one  correspondence  between the elements
of 
 $H^1(Z, {\mathcal{O}}_Z)$ 
 and  the above exact sequences. 
Further,  we have commutative diagram

\[
  \begin{array}{ccccccccc}
0\longrightarrow{\mathcal{O}}_Z &
\longrightarrow L^{\oplus n}&
\longrightarrow G &\longrightarrow 0\\
\quad \quad   \parallel &
\quad\Big\uparrow & 
\quad \Big\uparrow  & \\
0\longrightarrow{\mathcal{O}}_Z &\longrightarrow E & 
 \longrightarrow {\mathcal{O}}_Z &
\longrightarrow 0. 
\end{array}
\]

Since ${\Bbb{C}}\subset  H^1(Z, {\mathcal{O}}_Z )=
{\Bbb{C}}[x, \frac{1}{x}, \frac{1}{x-1}]$,
  $1\in H^1(Z, {\mathcal{O}}_Z )$. 
This nonzero element  $1$  corresponds to a  rank 2 vector bundle  $E$
such that  when restricted to 
every fibre $C_t$, it is the nonsplit extension of 
${\mathcal{O}}_{C_t}$  by ${\mathcal{O}}_{C_t}$.  
In fact, in the natural restriction map 
$$  H^1(Z, {\mathcal{O}}_Z )
\longrightarrow 
H^1(C_t, {\mathcal{O}}_{C_t}),
$$
 $1$  goes to  $1$.  A nonzero element of 
$H^1(C_t, {\mathcal{O}}_{C_t})$  determines a nonsplit 
extension of ${\mathcal{O}}_{C_t}$ by ${\mathcal{O}}_{C_t}$.
\begin{flushright}
 Q.E.D. 
\end{flushright}

\begin{lemma}There is a divisor $D$ on 
$X={\Bbb{P}}_Z(E)$ such that when restricted to $X_t={\Bbb{P}}_{C_t}(E_t)$,
$D|_{X_t}=D_t$ is the canonical section of $X_t$.
\end{lemma}
$Proof.$ By the above lemma,  we have surjective map from $E$ to 
${\mathcal{O}}_Z$. It corresponds to a section 
$\sigma : Z\longrightarrow  X.$ When restricted to $C_t$,
$\sigma|_{C_t}=\sigma_t : C_t\longrightarrow  X_t$
is the unique nonsplit extension of ${\mathcal{O}}_{C_t}$ by
${\mathcal{O}}_{C_t}$. 
\begin{flushright}
 Q.E.D. 
\end{flushright}

Let $Y=X-D$,
by the proof Theorem 1.3, we have $H^i(Y, \Omega^j_Y)=0$
for all $i>0$ and $j\geq 0$. 
By now we have constructed  a nonaffine, nonproduct 
example of threefold $Y$ with vanishing Hodge cohomology. 
This proves Theorem 1.6.

In the example, every fibre is  Stein and the base curve 
is Stein but we do not know whether the threefold is Stein. 
It is also interesting to construct  a threefold  with type (3) 
open surfaces as fibres.

 \end{document}